\documentclass{article}
\usepackage{amsmath}
\usepackage[latin1]{inputenc}
\usepackage{amssymb}
\usepackage{amscd}
\usepackage{natbib}
\usepackage{mathrsfs}
\usepackage[pdftex]{graphicx}
\usepackage{graphicx}
\usepackage[OT2, T1]{fontenc}
\usepackage[russian,english]{babel}
\usepackage[babel=true]{csquotes}
\usepackage{chngcntr}
\counterwithout{figure}{section}
\counterwithout{figure}{subsection}

\newcommand{\One}[0]{\ensuremath{\mathbf{1}}}

\newcommand{\R}[0]{\ensuremath{\mathbb{R}}}
\newcommand{\Z}[0]{\ensuremath{\mathbb{Z}}}

\DeclareMathOperator{\esssup}{ess\,sup}

\newtheorem{theo}{\sc{Theorem}}[section]
\newtheorem{lemm}[theo]{\sc{Lemma}}
\newtheorem{defi}[theo]{\sc{Definition}}

\newtheorem{cor}[theo]{\sc{Corollary}}
\newtheorem{rmq}[theo]{\sc{Remark}}

\begin{document}

\title{Decaying Turbulence in Generalised Burgers Equation}
\author{Alexandre Boritchev}
\date{\today}

\maketitle

\begin{center}
Department of Theoretical Physics (DPT)
\\
University of Geneva
\\
1211 GENEVA 4
\\
SWITZERLAND
\\
E-mail: alexandre.boritchev@gmail.com
\\
Telephone number: (+41) 22 37 96 312
\\
Fax number: (+41) 22 37 96 870
\end{center}

\textbf{Abstract.}\ We consider the generalised Burgers equation
\begin{equation} \nonumber
\frac{\partial u}{\partial t} + f'(u)\frac{\partial u}{\partial x} - \nu \frac{\partial^2 u}{\partial x^2}=0,\ t \geq 0,\ x \in S^1,
\end{equation}
where $f$ is strongly convex and $\nu$ is small and positive. We obtain sharp estimates for Sobolev norms of $u$ (upper and lower bounds differ only by a multiplicative constant). Then, we obtain sharp estimates for small-scale quantities which characterise the decaying Burgers turbulence, i.e. the dissipation length scale, the structure functions and the energy spectrum. The proof uses a quantitative version of an argument by Aurell, Frisch, Lutsko and Vergassola \cite{AFLV92}.
\\ \indent
Note that we are dealing with \textit{decaying}, as opposed to stationary turbulence. Thus, our estimates are not uniform in time. However, they hold on a time interval $[T_1, T_2]$, where $T_1$ and $T_2$ depend only on $f$ and the initial condition, and do not depend on the viscosity.
\\ \indent
These results give a rigorous explanation of the one-dimensional Burgers turbulence in the spirit of Kolmogorov's 1941 theory. In particular, we obtain two results which hold in the inertial range. On one hand, we explain the bifractal behaviour of the moments of increments, or structure functions. On the other hand, we obtain an energy spectrum of the form $k^{-2}$. These results remain valid in the inviscid limit.

\tableofcontents

\section{Introduction}

\subsection{Setting} \label{introset}
\begin{figure} 
\includegraphics[width=9cm, height=5cm]{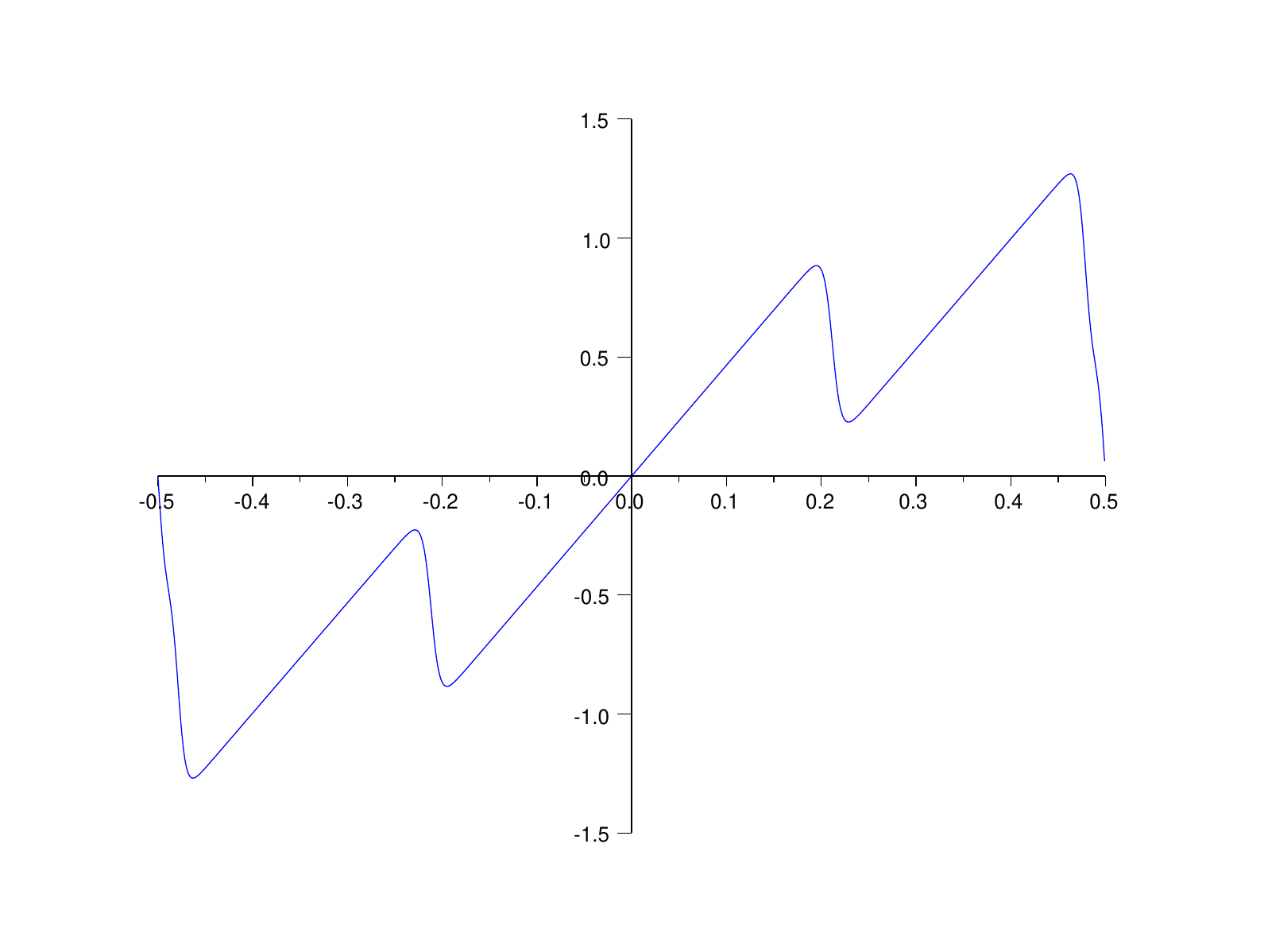}
\caption{\protect\label{N}  \enquote{Typical} solution of the Burgers equation}
\end{figure}

\smallskip
\indent
The Burgers equation
\begin{equation} \label{classBurgers}
\frac{\partial u}{\partial t} + u \frac{\partial u}{\partial x} - \nu \frac{\partial^2 u}{\partial x^2} = 0,
\end{equation}
where $\nu>0$ is a constant, appears in many fields of physics and other branches of science: see the reviews \cite{BF01,BK07} and references therein.
\\ \indent
The Burgers equation has been mentioned for the first time by Forsyth \cite{For06} and Bateman \cite{Bate15}, in 1906 and 1915 respectively. However, it only became well-known in the physical community around 1950, due to the work of the physicist whose name was given to it (see the monograph \cite{Bur74} and references therein). Burgers considered this equation as a toy model for hydrodynamics: indeed, the incompressible Navier-Stokes equations and (\ref{classBurgers}) have similar nonlinearities and dissipative terms, so this equation can be seen as the most natural one-dimensional model for Navier-Stokes.
\\ \indent
The equation (\ref{classBurgers}) can be transformed into the heat equation by the Cole-Hopf transformation \cite{Col51, Hop50}. However, this transformation will not be used in this paper for two different reasons. On one hand, the resulting representation of the solution is very singular as $\nu \rightarrow 0^+$, and interpreting this singularity rigorously is highly non-trivial. On the other hand, we want to be able to study the Burgers equation with $u\ \partial u /\partial x$ replaced by a more general nonlinearity; see (\ref{Burbegin})-(\ref{strconvex}).
\\ \indent
For $\nu \ll 1$, solutions of the Burgers equation display non-trivial small-scale behaviour, often referred to as decaying Burgers turbulence or \enquote{Burgulence} \cite{Bur74,Cho75,Kid79}. The language of the Kolmogorov 1941 theory \cite{Kol41a,Kol41b,Kol41c} is traditionally used to describe this behaviour. 
\\ \indent
For simplicity, from now on we consider the space-periodic setting, i.e. $x \in S^1=\R/\Z$. In this setting, the solutions of (\ref{classBurgers}) remain of order $1$ during a time of order $1$. On the other hand, for $t \rightarrow +\infty$ the solutions decay at least as $C t^{-1}$ in any Lebesgue space $L_p,\ 1 \leq p \leq +\infty$, uniformly in $\nu$ (cf. for instance \cite{Kru64}). Note that in the limit $\nu \rightarrow 0$, the diffusive effect due to the second derivative vanishes and this upper bound becomes sharp \cite[Theorem 11.7.3]{Daf10}. Thus, the solutions display smooth ramps and sharp cliffs \cite{BF01}. In the limit $\nu \rightarrow 0$, they have the $N$-wave behaviour, i.e. solutions are composed of waves similar to the Cyrillic capital letter \foreignlanguage{russian}{I} (the mirror image of N). In other words, at a fixed (large enough) time $t$ the solution $u(t,\cdot)$ alternates between negative jump discontinuities and smooth regions where the derivative is positive and of the order $1$ (see for instance \cite{Eva08}). This is a clear manifestation of the \text{small-scale intermittency} in space \cite{Fri95}. For $0<\nu \ll 1$ the solutions are still highly intermittent: there are zones where the derivative is small and positive, called {\it ramps},  and zones where the derivative is large in absolute value and negative, called {\it cliffs}.
\\ \indent
For a typical initial data $u_0$ (i.e. for $\max |u_0| \sim 1$ and $\max |(u_0)_x| \sim 1$) and for $t > 1/(\min (u_0)_x),\ t \sim 1$, it is numerically observed \cite{AFLV92} that a solution $u(t,\cdot)$ has the following features (cf. Figure~\ref{N}):
\begin{itemize}
\item{Amplitude of the solution:} $\sim 1$.
\\
\item{Number of cliffs per period:} $\sim 1$.
\\
\item{\enquote{Vertical drop} at a cliff:} $\sim -1$.
\\
\item{\enquote{Width} of a cliff:} $\sim \nu$.
\end{itemize}
\smallskip \indent
\smallskip
It is easy to verify that for the prototypical $N$-wave, i.e. for the $1$-periodic function equal to $x$ on $(-1/2,\ 1/2]$, the Fourier coefficients satisfy $|\hat{u}(k)| \sim k^{-1}$. Thus, it is natural to conjecture that for $\nu$ small and for a certain range of wave numbers $k$, the energy-type quantities $\frac{1}{2} |\hat{u}(k)|^2$ behave, in average, as $k^{-2}$ \cite{Cho75,FouFri83,Kid79,Kra68}.
\\ \indent
In the physical space, the natural analogues of the small-scale quantities  $\frac{1}{2} |\hat{u}(k)|^2$ are the structure functions
$$
S_p(\ell)=\int_{S^1}{|u(x+\ell)-u(x)|^p\ dx}.
$$
For $p \geq 0$, the description above implies that for $\nu \ll \ell \ll 1$, these quantities behave as $\ell^{\max(1,p)}$: in other words, we have a \textit{bifractal behaviour} \cite[Chapter 8]{Fri95}. Indeed, as observed in \cite{AFLV92}, there are three possibilities for the interval $[x,x+\ell]$ (below, $C$ denotes constants of order $1$, and we have to keep in mind that $\ell \ll 1$):
\smallskip
\begin{itemize}
\item{$[x,x+\ell]$ covers a large part of a "cliff".}
\\
Since the number of cliffs per period is of order $1$ and $\ell$ is larger than a cliff, the probability of this event is of order $\ell $. In this case:
$$
u(x+\ell)-u(x) \sim \underbrace{-C}_{"cliff"}+\underbrace{C \ell}_{"ramps"} \sim -C;\ |u(x+\ell)-u(x) |^p \sim C(p).
$$
\\
\item{$[x,x+\ell]$ covers a small part of a "cliff".}
\\
The contribution due to this possibility is negligible.
\\
\item{$[x,x+\ell]$ does not intersect a "cliff".}
\\
Since $\ell$ is smaller than the length of a ramp, the probability of this event is of order $1-C \ell \sim 1$. In this case:
$$
u(x+\ell)-u(x) \sim \underbrace{C \ell}_{ramp};\ |u(x+\ell)-u(x) |^p \sim C(p) \ell^p.
$$
\end{itemize}
\smallskip
\indent
Thus, for $\nu \ll \ell \ll 1$,
$$
S_p(\ell) \sim C(p) \ell+\ell^p \sim \left\lbrace \begin{aligned} & C(p) \ell^p,\ 0 \leq p \leq 1. \\ & C(p) \ell,\ p \geq 1. \end{aligned} \right.
$$

\subsection{Burgers equation and turbulence} \label{introturb}
\smallskip
 \indent
From now on, we consider the generalised one-dimensional space-periodic Burgers equation
\begin{equation} \label{Burbegin}
\frac{\partial u}{\partial t} + \frac{d f(u)}{d x} - \nu \frac{\partial^2 u}{\partial x^2} = 0,\quad x \in S^1=\R/\Z,
\end{equation}
where $f$ is $C^{\infty}$-smooth and strongly convex, i.e. $f$ satisfies the property
\begin{equation} \label{strconvex}
f''(y) \geq \sigma > 0,\quad y \in \R.
\end{equation}
The classical Burgers equation (\ref{classBurgers}) corresponds to $f(u)=u^2/2$. The physical arguments justifying the small-scale estimates which are given above still hold in that setting.
\\ \indent
For the sake of simplicity, we only consider solutions to (\ref{Burbegin})-(\ref{strconvex}) with zero space average for fixed $t$:
\begin{equation} \label{zero}
\int_{S^1}{u(t,x) dx}=0,\quad \forall t \geq 0.
\end{equation}
\medskip
\\ \indent
For the generalised Burgers equation, some \textit{upper} estimates for small-scale quantities have been obtained previously. Lemma~\ref{uxpos} of our paper is an analogue in the periodic setting of the one-sided Lipschitz estimate due to Oleinik, and the upper estimate for $S_1(\ell)$ follows from an estimate for the solution in the class of bounded variation functions $BV$. For references on these classical aspects of the theory of scalar conservation laws, see \cite{Daf10,Lax06,Ser99}. For some upper estimates for small-scale quantities, see \cite{Bar78,Kre88,Tad93}.
\\ \indent
Estimating small-scale quantities for nonlinear PDEs with small viscosity from above \textit{and from below} is motivated by the problem of turbulence. This research was initiated by Kuksin, who obtained estimates for a large class of equations (see \cite{Kuk97GAFA,Kuk99GAFA} and the references in \cite{Kuk99GAFA}).
\\ \indent
In the paper \cite{Bir01}, Biryuk obtained lower and upper estimates for the $L_2$-Sobolev norms of solutions to (\ref{Burbegin}). These estimates are sharp, in the sense that the lower and the upper bounds only differ by a multiplicative constant. Moreover, he obtained upper and lower estimates for the energy spectrum which enable him to give the correct value for the dissipation length scale. In \cite{BorK,BorW}, based on a better understanding of solutions for small values of $\nu$, we obtain sharp results for $L_p$-Sobolev norms, $p \in (1,\infty]$, and small-scale quantities. However, in both articles we add a rough in time and smooth in space random forcing term in the right-hand side of equation (\ref{Burbegin}) (a \enquote{kicked} and a white force, respectively). Thus, we change the nature of the equation: the energy injection due to the random forcing now balances the dissipation due to the second derivative. In other words, we study \textit{stochastic stationary} Burgulence, which is different from \textit{decaying} Burgulence. 
\\ \indent
Note that it is also possible to study (\ref{Burbegin}) in a deterministic stationary setting, which amounts to considering a deterministic additive random force. However, this is a delicate issue: indeed, for any initial condition $u_0$ we can build a \enquote{bad} time-independent random force equal to $f'(u_0)(u_0)_x-\nu (u_0)_{xx}$, corresponding to a stationary solution of (\ref{Burbegin}) which manifests no turbulent behaviour.
\\ \indent
Here, we prove sharp lower and upper estimates for the small-scale quantities, i.e. for the dissipation length scale, the structure functions and the energy spectrum, which characterise the decaying Burgulence. Thus, we improve significantly the results of \cite{Bir01}. To our best knowledge, this is the first such result for the \textit{deterministic} generalised Burgers equation. Moreover, we extend the results for the $L_2$-Sobolev norms obtained by Biryuk to the $L_p$-Sobolev norms, $p \in (1,\infty]$. The powers of $\nu,\ell,k$ involved in our estimates turn out to be the same as in the randomly forced case considered in \cite{BorK,BorW}. Note that our estimates hold in average on a time interval $[T_1,T_2]$, where both $T_1$ and $T_2$ do not depend on $\nu$. In other words, we consider a time range during which we have the transitory behaviour which is referred to as decaying Burgers turbulence \cite{BK07}. This time interval depends only on $f$ and, through the quantity $D$ (see \ref{D})), on $u_0$. In particular, it does not depend on $\nu$.
\\ \indent
A detailed overview of the results mentioned above is given in Section~\ref{results} (for the state of art) and in Section~\ref{resultspaper} (for the main results in this paper).
\\ \indent
Note that when studying the typical behaviour for solutions of nonrandom PDEs, one usually considers some averaging in the initial condition in order to avoid pathological initial data. Indeed, unlike for the stochastic case, now there is no random mechanism to get solutions out of \enquote{bad} regions of the phase space. Here, no such averaging is necessary. This is due to the particular structure of the deterministic Burgers equation: a non-zero initial condition $u_0$ is as \enquote{generic} as the ratio between the orders of $(u_0)_x$ and of $u_0$ itself. This ratio can be bounded from above using the quantity $D$:
\begin{equation} \label{D}
D=\max (|u_0|_{1}^{-1},\ | u_0 |_{1,\infty} )>1
\end{equation}
(see Subsection~\ref{Sob} for the meaning of the notation $|\cdot|_{m,p}$). Note that for $0 \leq m \leq 1$ and $1 \leq p \leq \infty$, we have:
\begin{equation} \label{Dp}
D^{-1} \leq |u_0|_{m,p} \leq D.
\end{equation}
The physical meaning of $D$ is that it gives a lower bound for the ratio between the amount of energy $\frac{1}{2} \int_{S^1}{u^2}$ initially contained in the system and its rate of dissipation $\nu \int_{S^1}{u_x^2}$.
\\ \indent
Now let us say a few words about similarities and differences between the Burgulence and real turbulence. It is clear that the geometric structures which are responsible for non-trivial small-scale behaviour are quite different for these two models: $N$-waves do not have the same properties as complex multi-scale structures such as vortex tubes observed in the real turbulence. However, because of the similarity in the form of the Burgers equation and the Navier-Stokes equations, the physical arguments justifying different theories of turbulence can be applied to the Burgulence. Indeed, both models exhibit an inertial nonlinearity of the form $u \cdot \nabla u$, and viscous dissipation which in the limit $\nu \rightarrow 0$ gives a dissipative anomaly \cite{BK07}. Hence, the Burgers equation is often used as a benchmark for the turbulence theories. It is also used as a benchmark for different numerical methods for the Navier-Stokes equations. For more information on both subjects, see \cite{BK07}.
\\ \indent
Now consider the generalised Burgers equation with a random regular in space and white in time forcing term $\eta$ such as in \cite{KuSh12}. Then the generalised Burgers equation with the natural scaling for this term (needed to counterbalance the energy dissipation due to the viscous term) is of the form:
$$
u_t+f'(u)u_x=\nu u_{xx}+\eta,
$$
i.e. the force does not depend on $\nu$ \cite{BorW}. This is similar to the conjectured behaviour for real turbulence, and contrasts with the situation for the 2D Navier-Stokes equations, where the corresponding term is of the form $\nu^{1/2} \eta$ \cite{KuSh12}. This justifies the study of the small-scale quantities for the randomly forced Burgers equation in the limit $\nu \rightarrow 0$ such as in \cite{BorW}. As it will be shown in Section~\ref{turb}, on a time scale which only depends on the initial condition and on the form of the nonlinearity $f'(u)u_x$, the small-scale quantities for the unforced Burgers equation also have a non-trivial behaviour as $\nu \rightarrow 0$, similar to the behaviour in the stochastic case. This is the main result of the paper. Up to now this question has only been adressed rigorously by Biryuk \cite{Bir01}, who obtained less sharp estimates. For more details on his results, see Section~\ref{results}.

\subsection{Plan of the paper}
\smallskip
 \indent
We introduce the notation and the setup in Section~\ref{nota}. In Section~\ref{results}, we give an overview of the state of art, before presenting the main results of our paper in Section~\ref{resultspaper}.
\\ \indent
In Section~\ref{Sobolev}, we begin by recalling an upper estimate for the quantity $\partial u/\partial x$.  This result allows us to obtain upper bounds, as well as time-averaged lower bounds, for the Sobolev norms $|u|_{m,p}$. These bounds depend only on $f$ and on the quanity $D$ defined by (\ref{D}).
\\ \indent
In Section~\ref{turb} we give sharp upper and lower bounds for the dissipation length scale, the structure functions and the energy spectrum for the flow $u(t,x)$, which hold uniformly for $\nu \leq \nu_0$, and we analyse the meaning of these results in terms of the theory of turbulence. These bounds and the constant $\nu_0>0$ only depend on $f$ and on $D$.
\\ \indent
In Section~\ref{inviscid} we consider the inviscid limit $\nu=0$.

\section{Notation and setup} \label{nota}
\smallskip
 \indent
\textbf{Agreement}: In the whole paper, all functions that we consider are real-valued and the space variable $x$ belongs to $S^1=\R/\Z$.

\subsection{Sobolev spaces} \label{Sob}
\smallskip
 \indent
Consider a zero mean value integrable function $v$ on $S^1$. 
For $p \in [1,\infty)$, we denote its $L_p$ norm
$$
\Bigg( \int_{S^1}{|v|^p} \Bigg)^{1/p}
$$
by $\left|v\right|_p$. The $L_{\infty}$ norm is by definition
$$
\left|v\right|_{\infty}=\esssup_{x \in S^1} |v(x)|.
$$
The $L_2$ norm is denoted by  $\left|v\right|$, and $\left\langle \cdot,\cdot\right\rangle$ stands for the $L_2$ scalar product. From now on $L_p,\ p \in [1,\infty],$ denotes the space of zero mean value functions in $L_p(S^1)$. Similarly, $C^{\infty}$ is the space of $C^{\infty}$-smooth zero mean value functions on $S^1$.
\\ \indent
For a nonnegative integer $m$ and $p \in [1,\infty]$, $W^{m,p}$ stands for the Sobolev space of zero mean value functions $v$ on $S^1$ with finite norm
\begin{equation} \nonumber
\left|v\right|_{m,p}=\left|\frac{d^m v}{dx^m}\right|_p.
\end{equation}
In particular, $W^{0,p}=L_p$ for $p \in [1,\infty]$. For $p=2$, we denote $W^{m,2}$ by $H^m$, and abbreviate the corresponding norm as $\left\|v\right\|_m$. 
\\ \indent
Note that since the length of $S^1$ is $1$ and the mean value of $v$ vanishes, we have:
$$
|v|_1 \leq |v|_{\infty} \leq |v|_{1,1} \leq |v|_{1,\infty} \leq \dots \leq |v|_{m,1} \leq |v|_{m,\infty} \leq \dots
$$
We recall a version of the classical Gagliardo-Nirenberg inequality: cf. \cite[Appendix]{DG95}.
\begin{lemm} \label{GN}
For a smooth zero mean value function $v$ on $S^1$,
$$
\left|v\right|_{\beta,r} \leq C \left|v\right|^{\theta}_{m,p} \left|v\right|^{1-\theta}_{q},
$$
where $m>\beta$, and $r$ is determined by
$$
\frac{1}{r}=\beta-\theta \Big( m-\frac{1}{p} \Big)+(1-\theta)\frac{1}{q},
$$
under the assumption $\theta=\beta/m$ if $p=1$ or $p=\infty$, and $\beta/m \leq \theta < 1$ otherwise. The constant $C$ depends on $m,p,q,\beta,\theta$.
\end{lemm}
\indent
Subindices $t$ and $x$, which can be repeated, denote partial differentiation with respect to the corresponding variables. We denote by $v^{(m)}$ the $m$-th derivative of $v$ in the variable $x$. The function $v(t,\cdot)$ is abbreviated as $v(t)$.

\subsection{Notation} 
\smallskip
 \indent
In this paper, we study asymptotical properties of solutions to (\ref{Burbegin}) for small values of $\nu$, i.e. we suppose that
$$
0 < \nu \ll 1.
$$
We assume that $f$ is infinitely differentiable and satisfies (\ref{strconvex}). We recall that we restrict ourselves to the zero space average case, i.e. the initial condition $u_0:=u(0)$ satisfies (\ref{zero}). Consequently, $u(t)$ satisfies (\ref{zero}) for all $t$. Furthermore, we assume that $u_0 \in C^{\infty}$. We also assume that we are not in the case $u_0 \equiv 0$, corresponding to the trivial solution $u(t,x) \equiv 0$. This ensures that the quantity $D$ (see (\ref{D})) is well-defined.
\\ \indent
For the existence, uniqueness and smoothness of solutions to (\ref{Burbegin}), see for instance \cite{KrLo89}.
\smallskip
\\ \indent
\textbf{Agreements:}\ From now on, all constants denoted by $C$ with sub- or superindexes are positive. Unless otherwise stated, they depend only on $f$ and on $D$. By $C(a_1,\dots,a_k)$ we denote constants which also depend on parameters $a_1,\dots,a_k$. By $X \overset{a_1,\dots,a_k}{\lesssim} Y$ we mean that $X \leq C(a_1,\dots,a_k) Y$. The notation $X \overset{a_1,\dots,a_k}{\sim} Y$ stands for
$$
Y \overset{a_1,\dots,a_k}{\lesssim} X \overset{a_1,\dots,a_k}{\lesssim} Y.
$$
In particular, $X \lesssim Y$ and $X \sim Y$ mean that $X \leq C Y$ and $C^{-1} Y \leq X \leq C Y$, respectively.
\\ \indent
All constants are independent of the viscosity $\nu$. We denote by $u=u(t,x)$ a solution of (\ref{Burbegin}) for an initial condition $u_0$. A relation where the admissible values of $t$ (respectively, $x$) are not specified is assumed to hold for all $t \geq 0$ or $t>0$, depending on the context (respectively, all $x \in S^1$).
\\ \indent
The brackets $\lbrace \cdot \rbrace$ stand for the averaging in time over an interval $[T_1,T_2]$, where $T_1,T_2$ only depend on $f$ and on $D$ (see (\ref{T1T2}) for their definition.)
\\ \indent
For $m \geq 0$, $p \in [1,\infty]$, $\gamma(m,p)$ is by definition the quantity $\max(0,m-1/p)$.
\\ \indent
We use the notation $g^{-}=\max(-g,0)$ and $g^{+}=\max(g,0)$.

\subsection{Notation in Section~\ref{turb}} \label{notaturb}
\smallskip
\indent
In that section, we study analogues of quantities which are important for hydrodynamical turbulence. We consider quantities in physical space (structure functions) as well as in Fourier space (energy spectrum). We assume that $\nu \leq \nu_0$. The value of $\nu_0>0$ will be chosen in (\ref{nu0eq}).
\\
\smallskip
\begin{figure} 
\includegraphics[height=2cm]{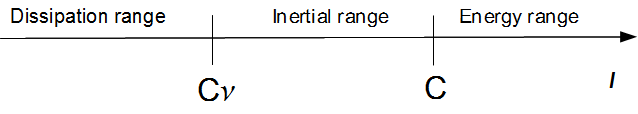}
\caption{\protect\label{figscales} Scales for the 1D Burgers solutions}
\end{figure}
\indent
We define the non-empty and non-intersecting intervals
\begin{equation} \nonumber
J_1=(0,\ C_1 \nu];\ J_2=(C_1 \nu,\ C_2];\ J_3=(C_2,\ 1]
\end{equation}
(see Figure~\ref{figscales}) corresponding to the \textit{dissipation range}, the \textit{inertial range} and the \textit{energy range} from the Kolmogorov 1941 theory of turbulence, respectively \cite{Fri95}. In particular, the upper bound $C_1 \nu$ of the dissipation range is the \textit{dissipation length scale}. The positive constants $C_1$ and $C_2$ will be chosen in (\ref{K}-\ref{nu0eq}) in such a manner that $C_1 \nu_0<C_2<1$, which ensures that the intervals $J_i$ are non-empty and non-intersecting.
\\ \indent
By Theorem \ref{avoir} we obtain that $\lbrace |u|^2 \rbrace \sim 1$ (see above for the meaning of the notation $\lbrace \cdot \rbrace$). On the other hand, by (\ref{W11}) (after integration by parts) we get:
\begin{align} \label{intparts}
\lbrace |\hat{u}(n)|^2 \rbrace &= (2 \pi n)^{-2} \Big\{ \Big| \int_{S^1}{e^{2 \pi i n x} u_x(x)} \Big|^2 \Big\} \leq (2 \pi n)^{-2} \lbrace |u|^2_{1,1} \rbrace \leq C n^{-2},
\end{align}
and $C_1$ and $C_2$ can be made as small as we wish (see (\ref{nu0ineq})). Consequently, the proportion of the sum $\lbrace\sum |\hat{u}(n)|^2 \rbrace$ contained in the Fourier modes corresponding to $J_3$ can be made as large as desired. For instance, we may assume that
$$
\Bigg\{ \sum_{|n| < C_2^{-1}}{|\hat{u}(n)|^2} \Bigg\} \geq \frac{99}{100} \Bigg\{ \sum_{n \in \Z}{|\hat{u}(n)|^2} \Bigg\}.
$$
The quantities $S_{p}(\ell)$ denote the averaged moments of the increments in space for the flow $u(t,x)$:
$$
S_{p}(\ell)= \Bigg\{ \int_{S^1}{|u(t,x+\ell)-u(t,x)|^p dx} \Bigg\},\ p \geq 0,\ 0< \ell \leq 1.
$$
The quantity $S_{p}(\ell)$ is the \textit{structure function} of $p$-th order. The flatness, which measures spatial intermittency \cite{Fri95}, is defined by:
\begin{equation} \label{flatness}
F(\ell)=S_4(\ell)/S_2^2(\ell).
\end{equation}
Finally, for $k \geq 1$, we define the (layer-averaged) energy spectrum by
\begin{equation} \label{spectrum}
E(k)=\Bigg\{ \frac{\sum_{|n| \in [M^{-1}k,Mk]}{|\hat{u}(n)|^2}}{\sum_{|n| \in [M^{-1}k,Mk]}{1}} \Bigg\},
\end{equation}
where $M \geq 1$ is a constant which will be specified later (see the proof of Theorem~\ref{spectrinert}).

\section{State of the art} \label{results}
\smallskip
\indent
We recall that $u=u(t,x)$ denotes a solution of (\ref{Burbegin}) for an initial condition $u_0$. All constants are independent of the viscosity $\nu$ (i.e., dependance on $\nu$ is always explicitly stated in the estimates). A relation where the admissible values of $t$ (respectively, $x$) are not specified is assumed to hold for all $t \geq 0$ or $t>0$, depending on the context (respectively, all $x \in S^1$). For more information on the notation, see Section~\ref{nota}.
\\ \indent
The estimate
\begin{equation}\label{E}
u_x(t,x) \leq (\sigma t)^{-1},\quad t>0,
\end{equation}
where $\sigma$ is the constant in the assumption (\ref{strconvex}), is a reformulation of Oleinik's $E$-condition \cite{Ole57}. This result immediately implies an upper bound for the first structure function $S_1(\ell)$. Indeed, since the space average of $u(t,\cdot)$ vanishes identically for all $t$, we have:
\begin{align} \nonumber
& \int_{S^1}{|u(t,x+\ell)-u(t,x)|} 
\\ \nonumber
&= \int_{S^1}{(u(t,x+\ell)-u(t,x))^{+}}+\int_{S^1}{(u(t,x+\ell)-u(t,x))^{-}}
\\ \nonumber
&= 2 \int_{S^1}{(u(t,x+\ell)-u(t,x))^{+}} \leq \frac{2}{\sigma t} \ell.
\end{align}
Moreover, integration by parts gives us the follwing upper estimate for the spectrum:
$$
\lbrace |\hat{u}(n)|^2 \rbrace \leq C (\sigma t n)^{-2}
$$
(see for instance \cite{Bar78}). In a similar setting, exponential upper estimates for the spectrum in the dissipation range have also been obtained; see \cite{Kre88}. See also \cite{Tad93} for upper estimates in a slightly different (hyperviscous) setting.
\bigskip
\\ \indent
In \cite{Bir01}, Biryuk begins by proving upper as well as \textit{lower} estimates for the $H^m$-Sobolev norms of $u$:
\begin{equation} \label{BirSob}
B^{-1} \nu^{-(2m-1)} \leq \frac{1}{T} \int_0^T{\Vert u \Vert_m^2} \leq B \nu^{-(2m-1)}, \quad 0 < \nu \leq \nu_0.
\end{equation}
Here,  the strictly positive quantities $\nu_0$ and $T$ depend on $f$ and $m$ as well as on the Sobolev norms of the initial condition $u_0$. The letter $B$ denotes different strictly positive quantities which also depend on these parameters. Since these estimates hold only for a fixed value of $T$, there is no contradiction with the decay in $C t^{-1}$ of the solutions as $t \rightarrow +\infty$.
\\ \indent
Let us denote by $E_{s,\theta}$ the averaged energy spectrum:
$$
E_{s,\theta}=\frac{1}{T} \int_0^T{ \frac{\sum_{|n| \in [\nu^{-s+\theta}, \nu^{-s-\theta}]}{|\hat{u}(n)|^2}}{\sum_{|n| \in [\nu^{-s+\theta}, \nu^{-s-\theta}]}{1}} },\quad s,\theta>0.
$$
Using the inequalities (\ref{E})-(\ref{BirSob}), Biryuk obtains upper and lower estimates for the spectrum of the solutions, which hold for $0 < \nu \leq \nu_0$:
\begin{align} \label{Birspecdiss}
& E_{s,\theta} \leq B \nu^{m},\quad m>0,\ s>1+\theta.
\\ \label{Birspecupp}
& E_{s,\theta} \leq B \nu^{2(s-\theta)}, \quad s>\theta.
\\ \label{Birspeclow}
& E_{1,\theta} \geq B \nu^{2+2 \theta}.
\end{align}
The quantities $\nu_0$ and $T$, as well as the different strictly positive quantities denoted by $B$, depend on $f$ and on the Sobolev norms of $u_0$, as well as on $m$, $s$, $\theta$.
\\ \indent
Note that Biryuk's results for the Sobolev norms are sharp, in the sense that in the lower and upper estimates in (\ref{BirSob}), $\nu$ is raised to the same power. Using the same terminology, his results (\ref{Birspecupp}-\ref{Birspeclow}) can be described as "almost sharp" for $s=1$, since they give almost the same lower and upper bounds for $E_{1,\theta}$ with $0 < \theta \ll 1$ (up to a multiplicative constant and $\nu$ raised to a very small power).
\\ \indent
Biryuk's spectral estimates may be interpreted in the spirit of Kolmogorov's theory of turbulence. Indeed, relation (\ref{Birspecdiss}) implies that the energy spectrum of the $k$-th Fourier mode averaged around $k=K$, where $K \gg \nu^{-1}$, decays faster than any negative degree of $K$. This suggests that for $K \gg \nu^{-1}$ we are in the dissipation range, where the energy $E_k$ decays fast. On the other hand, relations (\ref{Birspecupp}) and  (\ref{Birspeclow}) yield that the energy $E_k$, averaged around $k=\nu^{-1}$, behaves as $k^{-2}$, which gives a Kolmogorov-type power law \cite{Fri95}. This suggests a dissipation length scale of the order $\nu$.

\section{Main results} \label  {resultspaper}
\smallskip
 \indent
In our paper, in Section~\ref{Sobolev}, we prove sharp upper and lower bounds for almost all Sobolev norms of $u$, generalising the estimates (\ref{BirSob}). These results for Sobolev norms of solutions are summed up in Theorem~\ref{avoir}. Namely, for $m \in \lbrace 0,1 \rbrace$ and $p \in [1,\infty]$ or for $m \geq 2$ and $p \in (1,\infty]$ we have:
\begin{equation} \label{avoirresults}
\Big( \lbrace {\left|u(t)\right|_{m,p}^{\alpha}} \rbrace \Big)^{1/\alpha} \overset{m,p,\alpha}{\sim} \nu^{-\gamma},\quad \alpha>0.
\end{equation}
We recall that by definition, $\gamma(m,p)=\max(0,m-1/p)$, and the brackets $\lbrace \cdot \rbrace$ stand for the averaging in time over an interval $[T_1,T_2]$ ($T_1, T_2$ only depend on $f$ and, through $D$, on $u_0$: see (\ref{T1T2})). For more information on the notation, see Section~\ref{nota}.
\\ \indent
In Section~\ref{turb} we obtain sharp estimates for analogues of quantities characterising hydrodynamical turbulence. In what follows, we assume that $\nu \in (0,\nu_0]$, where $\nu_0 \in (0,1]$ depends only on $f$ and on $D$.
\\ \indent
First, as a consequence of (\ref{E}) and (\ref{avoirresults}), in Theorem~\ref{avoir2} we prove that for $\ell \in J_1$:
$$
\quad \ \ \ S_{p}(\ell) \overset{p}{\sim} \left\lbrace \begin{aligned} & \ell^{p},\ 0 \leq p \leq 1. \\ & \ell^{p} \nu^{-(p-1)},\ p \geq 1, \end{aligned} \right.$$
and for $\ell \in J_2$:
$$
S_{p}(\ell) \overset{p}{\sim} \left\lbrace \begin{aligned} & \ell^{p},\ 0 \leq p \leq 1. \\ & \ell,\ p \geq 1. \end{aligned} \right.
$$
Consequently, for $\ell \in J_2$ the flatness satisfies the estimate:
$$
F(\ell)=S_4(\ell)/S_2^2(\ell) \sim \ell^{-1}.
$$
Thus, $u$ is highly intermittent in the inertial range. This intermittency is in good agreement with the physical heuristics presented in Subsection~\ref{introset}, due to the particular structure of the solution, where the excited zones correspond to the cliffs. Cf. \cite{Fri95} for a discussion of the intermittency for hydrodynamical turbulent flows.
\\ \indent
Finally, as a relatively simple consequence of our estimates for the structure function $S_2(\ell)$, we get estimates for the spectral asymptotics of the decaying Burgulence. On one hand, as a consequence of Theorem~\ref{avoir}, for $m \geq 1$ we get:
$$
\lbrace |\hat{u}(k)|^2 \rbrace \overset{m}{\lesssim} k^{-2m} {\Vert u \Vert_m^2 } \overset{m}{\lesssim} (k \nu)^{-2m} \nu.
$$
In particular, $\lbrace |\hat{u}(k)|^2 \rbrace$ decreases at a faster-than-algebraic rate for $|k| \succeq \nu^{-1}$. On the other hand, by Theorem~\ref{spectrinert}, for $k$ such that $k^{-1} \in J_2$ the energy spectrum $E(k)$ satisfies:
$$
E(k) \sim k^{-2},
$$
where the quantity $M \geq 1$ in the definition of $E(k)$ depends only on $f$ and on $D$. This result significantly improves Biryuk's spectral estimates, since it characterises exactly the spectral behaviour in the whole inertial range.
\\ \indent
Note that our estimates hold for quantities averaged on a time interval $[T_1,T_2],\ T_2>T_1>0$, and not on an interval $[0,T]$ as in Biryuk's paper. This allows us to obtain estimates which depend on the initial condition only through the single parameter $D$. Moreover, as in Biryuk's paper, this time interval does not depend on the viscosity coefficient $\nu$.
\\ \indent
As we mentioned in Section~\ref{introturb}, upper estimates for $S_p(\ell)$ follow from known results about the Burgers equation. Sharp lower estimates were not known before our work.
\\ \indent
Finally, in Section~\ref{inviscid} we note that our estimates for the small-scale quantities still hold in the inviscid limit $\nu \rightarrow 0$, up to some natural modifications.

\section{Estimates for Sobolev norms} \label{Sobolev}
\smallskip
 \indent
We recall that $u=u(t,x)$ denotes a solution of (\ref{Burbegin}) for an initial condition $u_0$. All constants are independent of the viscosity $\nu$. A relation where the admissible values of $t$ (respectively, $x$) are not specified is assumed to hold for all $t \geq 0$ or $t>0$, depending on the context (respectively, all $x \in S^1$). For more information on the notation, see Section~\ref{nota}.
\\ \indent
We begin by recalling a key upper estimate for $u_x$.

\begin{lemm} \label{uxpos}
We have:
$$
u_x(t,x) \leq \min (D,\sigma^{-1} t^{-1} ).
$$
\end{lemm}

\textbf{Proof.}
Differentiating the equation (\ref{Burbegin}) once in space we get
$$
(u_x)_t+f''(u)u_x^2+f'(u)(u_{x})_{x}=\nu (u_x)_{xx}.
$$
Now consider a point $(t_1,x_1)$ where $u_x$ reaches its maximum on the cylinder $S=[0,t] \times S^1$. Suppose that $t_1>0$ and that this maximum is nonnegative. At such a point, Taylor's formula implies that we would have $(u_x)_t \geq 0$, $(u_{x})_{x}=0$ and $(u_x)_{xx} \leq 0$. Consequently, since by (\ref{strconvex}) $f''(u) \geq \sigma$, we get $f''(u)u_x^2 \leq 0$, which is impossible. Thus $u_x$ can only reach a nonnegative maximum on $S$ for $t_1=0$. In other words, since $(u_0)_x$ has zero mean value, we have:
$$
u_x(t,x) \leq \max_{x \in S^1} {(u_0)_x(x)} \leq D.
$$
\\ \indent
The inequality
$$
u_x(t,x) \leq \sigma^{-1} t^{-1}
$$
is proved in by a similar maximum principle argument applied to the function $tu_x$: cf. \cite{Kru64}.\ $\square$
\smallskip
\\ \indent
Since the space averages of $u(t)$ and $u_x(t)$ vanish, we get the following upper estimates:
\begin{align} \label{Lpupper}
&\left|u(t)\right|_{p} \leq \left|u(t)\right|_{\infty} \leq \int_{S^1}{u_x^{+}(t)} \leq \min (D, \sigma^{-1} t^{-1} ),\quad 1 \leq p \leq +\infty,
\\ \label{W11}
&\left|u(t)\right|_{1,1}=\int_{S^1}{u_x^{+}(t)}+\int_{S^1}{u_x^{-}(t)}=2 \int_{S^1}{u_x^{+}(t)} \leq 2\min (D, \sigma^{-1} t^{-1}).
\end{align}

Now we recall a standard estimate for the nonlinearity $\left\langle v^{(m)}, (f(v))^{(m+1)}\right\rangle$. For its proof, we refer to \cite{BorW}.

\begin{lemm} \label{lmubuinfty}
For $v \in C^{\infty}$ such that $\left|v\right|_{\infty} \leq A$, we have:
$$
\left| \left\langle  v^{(m)}, (f(v))^{(m+1)} \right\rangle\right| \leq \tilde{C} \left\|v\right\|_m \left\|v\right\|_{m+1},\quad m \geq 1,
$$
where $\tilde{C}$ depends only on $m$, $A$ and $\left|f\right|_{C^m([-A,A])}$.
\end{lemm}
\indent
The following result shows that there is a strong nonlinear damping which prevents the successive derivatives of $u$ from becoming too large.

\begin{lemm} \label{uppermaux}
We have
$$
\left\|u(t)\right\|^{2}_1 \lesssim \nu^{-1}.
$$
On the other hand, for $m \geq 2$,
$$
\left\|u(t)\right\|^{2}_m \overset{m}{\lesssim} \max (\nu^{-(2m-1)}, t^{-(2m-1)}).
$$
\end{lemm}

\textbf{Proof.} Fix $m \geq 1$. Denote
$$
x(t)=\left\|u(t)\right\|^{2}_m.
$$
We claim that the following implication holds:
\begin{align} \label{decrm}
&x(t) \geq C' \nu^{-(2m-1)} \Longrightarrow \frac{d}{dt} x(t) \leq -(2m-1) x(t)^{2m/(2m-1)},
\end{align}
where $C'$ is a fixed positive number, chosen later. Below, all constants denoted by $C$ do not depend on $C'$.
\\ \indent
Indeed, assume that $ x(t) \geq C' \nu^{-(2m-1)}.$ Integrating by parts in space and using (\ref{Lpupper}) ($p=\infty$) and Lemma~\ref{lmubuinfty}, we get the following energy dissipation relation:
\begin{align} \nonumber
\frac{d}{dt}  x(t) &= - 2 \nu  \left\|u(t)\right\|_{m+1}^2-2\left\langle u^{(m)}(t), (f(u(t)))^{(m+1)}\right\rangle
\\ \indent
&\leq - 2 \nu  \left\|u(t)\right\|_{m+1}^2 + C \left\|u(t)\right\|_m \left\|u(t)\right\|_{m+1}.
\end{align}
Applying Lemma~\ref{GN} to $u_x$ and then using (\ref{W11}), we get:
\begin{align} \nonumber
\left\| u(t)\right\|_m &\leq C \left\| u(t)\right\|_{m+1}^{(2m-1)/(2m+1)} \left| u(t)\right|_{1,1}^{2/(2m+1)} 
\\ \label{GN11m}
&\leq C \left\| u(t)\right\|_{m+1}^{(2m-1)/(2m+1)}.
\end{align}
Thus, we have the relation
\begin{align} \label{int1}
\frac{d}{dt} x(t) \leq & (- 2 \nu \left\|u(t)\right\|_{m+1}^{2/(2m+1)} + C ) \left\|u(t)\right\|_{m+1}^{4m/(2m+1)}.
\end{align}
The inequality (\ref{GN11m}) yields
\begin{equation} \label{int2}
\left\| u(t)\right\|_{m+1}^{2/(2m+1)} \geq C x(t)^{1/(2m-1)},
\end{equation}
and then since by assumption $ x(t) \geq C' \nu^{-(2m-1)}$ we get:
\begin{align} \label{int3}
 \left\|u(t)\right\|_{m+1}^{2/(2m+1)} & \geq C C'^{1/(2m-1)} \nu^{-1}.
\end{align}
Combining the inequalities (\ref{int1}-\ref{int3}), for $C'$ large enough we get:
\begin{align} \nonumber
\frac{d}{dt}  x(t) &\leq (-C C'^{1/(2m-1)} + C ) x(t)^{2m/(2m-1)}.
\end{align}
Thus we can choose $C'$ in such a way that the implication (\ref{decrm}) holds.
\\ \indent
For $m=1$, (\ref{Dp}) and (\ref{decrm}) immediately yield that
$$
x(t) \leq \max(C'\nu^{-1},D^2) \leq \max(C',D^2) \nu^{-1},\ t \geq 0.
$$
Now consider the case $m \geq 2$. We claim that
\begin{equation} \label{decrmcor}
x(t) \leq \max ( C' \nu^{-(2m-1)}, t^{-(2m-1)} ).
\end{equation}
Indeed, if $x(s) \leq C' \nu^{-(2m-1)}$ for some $s \in \left[0,t\right]$, then the assertion (\ref{decrm}) ensures that $x(s)$ remains below this threshold up to time $t$.
\\ \indent
Now, assume that $x(s) > C' \nu^{-(2m-1)}$ for all $s \in \left[0,t\right]$. Denote
$$
\tilde{x}(s)=(x(s))^{-1/(2m-1)},\ s \in \left[0,t\right].
$$
By (\ref{decrm}) we get $d\tilde{x}(s)/ds \geq 1$. Therefore $\tilde{x}(t) \geq t$ and $x(t) \leq t^{-(2m-1)}$. Thus in this case, the inequality (\ref{decrmcor}) still holds. This proves the lemma's assertion. $\square$

\begin{lemm} \label{upperwmp}
For $m \in \lbrace 0,1 \rbrace$ and $p \in [1,\infty]$, or for $m \geq 2$ and $p \in (1,\infty]$ we have:
\begin{align} \nonumber
\left|u(t)\right|_{m,p} & \overset{m,p}{\lesssim} \max(\nu^{-\gamma},t^{-\gamma}).
\end{align}
\end{lemm}

\textbf{Proof.}
For $m \geq 1$ and $p \in [2,\infty]$, we interpolate $\left|u(t)\right|_{m,p}$ between $\left\|u(t)\right\|_{m}$ and $\left\|u(t)\right\|_{m+1}$. By Lemma~\ref{GN} applied to $u^{(m)}(t)$, we have:
$$
\left|u(t)\right|_{m,p} \overset{p}{\lesssim} \left\|u(t)\right\|_{m}^{1-\theta} \left\|u(t)\right\|_{m+1}^{\theta},\ \theta=\frac{1}{2}-\frac{1}{p}.
$$
Then we use Lemma~\ref{uppermaux} and H{\"o}lder's inequality to complete the proof.
\\ \indent
We use the same method to prove the case $m=1,\ p \in [1,2]$, combining (\ref{W11}) and Lemma~\ref{uppermaux}. We also proceed similarly for $m \geq 2,\ p \in (1,2)$, combining (\ref{W11}) and an upper estimate for $\Vert u(t) \Vert_{M,p}^{\alpha}$ for a large value of $M$ and some $p \geq 2$.
\\ \indent
Finally, the case $m=0$ follows from (\ref{Lpupper}). $\square$
\medskip
\\ \indent
Unfortunately, the proof of Lemma~\ref{upperwmp} cannot be adapted to the case $m \geq 2$ and $p=1$. Indeed, Lemma~\ref{GN} only allows us to estimate a $W^{m,1}$ norm from above by other $W^{m,1}$ norms: we can only get that
$$
|u(t)|_{m,1} \overset{m,n,k}{\lesssim} |u(t)|_{n,1}^{(m-k)/(n-k)} |u(t)|_{k,1}^{(n-m)/(n-k)},\ 0 \leq k < m < n,
$$
and thus the upper estimates obtained above cannot be used. However, we have:
$$
|u(t)|_{m,1} \leq |u(t)|_{m,1+\beta}
$$
for any $\beta>0$. Consequently, the lemma's statement holds for $m \geq 2$ and $p=1$, with $\gamma$ replaced by $\gamma+\lambda$, and $\overset{m,p}{\lesssim}$ replaced by $\overset{m,p,\lambda}{\lesssim}$, for any $\lambda>0$.
\medskip
\\ \indent
Now we define
\begin{equation} \label{T1T2}
T_1=\frac{1}{4}D^{-2} \tilde{C}^{-1};\quad T_2=\max \Big( \frac{3}{2} T_1,\quad 2D\sigma^{-1} \Big),
\end{equation}
where $\tilde{C}$ is a constant such that for all $t$, $\left\|u(t)\right\|_1^2 \leq \tilde{C} \nu^{-1}$ (cf. Lemma~\ref{uppermaux}). Note that $T_1$ and $T_2$ do not depend on the viscosity coefficient $\nu$. 
\\ \indent
From now on, for any function $A(t)$, $\lbrace A(t) \rbrace$ is by definition the time average
$$
\lbrace A(t) \rbrace=\frac{1}{T_2-T_1} \int_{T_1}^{T_2}{A(t)}.
$$
\\ \indent
The first quantity that we estimate from below is $\lbrace |u(t)|_p^2 \rbrace,\ p \in [1,\infty]$.

\begin{lemm} \label{finitetimep}
For $p \in [1,\infty]$, we have:
$$
\lbrace |u(t)|_p^2 \rbrace \gtrsim 1.
$$
\end{lemm}

\textbf{Proof.} It suffices to prove the lemma's statement for $p=1$. But this case follows from the case $p=2$. Indeed, by H{\"o}lder's inequality and (\ref{Lpupper}) we have:
$$
\lbrace |u(t)|_1^2 \rbrace \geq \lbrace |u(t)|_{\infty}^{-2} |u(t)|^{4} \rbrace \gtrsim \lbrace |u(t)|^{4} \rbrace \geq \lbrace |u(t)|^{2} \rbrace^2.
$$
\indent
Integrating by parts in space, we get the dissipation identity
\begin{align} \label{dissip}
\frac{d}{dt} \left| u(t)\right|^2 &= \int_{S^1}{(-2 uf'(u)u_x+2 \nu u u_{xx})}=-2 \nu \left\| u(t)\right\|_{1}^{2}.
\end{align}
Thus, integrating in time and using (\ref{D}) and Lemma~\ref{uppermaux}, we obtain that for $t \in [T_1,3T_1/2]$ we have the following uniform lower bound:
\begin{align} \label{dissipaux}
|u(t)|^2 &= |u_0|^2 - 2 \nu \int_{0}^{t}{\left\|u(t)\right\|_1^2} \geq D^{-2} - 3 T_1 \tilde{C} \geq D^{-2}/4.
\end{align}
Thus,
$$
\lbrace |u(t)|^2 \rbrace \geq  \frac{1}{T_2-T_1} \int_{T_1}^{3T_1/2}{|u(t)|^2} \geq  \frac{D^{-2} T_1}{8(T_2-T_1)}.\ \square
$$
\medskip \\ \indent
Now we prove a key estimate for $\lbrace \left\|u(t)\right\|_1^2 \rbrace$.

\begin{lemm} \label{finitetime}
We have
$$
\lbrace \left\|u(t)\right\|_1^2 \rbrace \gtrsim \nu^{-1}.
$$
\end{lemm}

\textbf{Proof.} Integrating (\ref{dissip}) in time in the same way as in (\ref{dissipaux}), we prove that $|u(T_1)|^2 \geq D^{-2}/2$. Thus, using (\ref{Lpupper}) ($p=2$) we get:
\begin{align} \nonumber
\lbrace \left\|u(t)\right\|_1^2 \rbrace &= \frac{1}{2 \nu (T_2-T_1)} (|u(T_1)|^2-|u(T_2)|^2)
\\ \nonumber
& \geq \frac{1}{2 \nu (T_2-T_1)} \Big(\frac{1}{2}D^{-2}-\sigma^{-2} T_2^{-2} \Big)
\\ \nonumber
& \geq \frac{D^{-2}}{8 (T_2-T_1) } \nu^{-1},
\end{align}
which proves the lemma's assertion.\ $\square$
\medskip \\ \indent
This time-averaged lower bound yields similar bounds for other Sobolev norms.

\begin{lemm} \label{finalexp}
For $m \geq 1$,
$$
\lbrace \left\|u(t)\right\|_m^2 \rbrace \overset{m}{\gtrsim} \nu^{-(2m-1)}.
$$
\end{lemm} 

\textbf{Proof.}
Since the case $m=1$ has been treated in the previous lemma, we may assume that $m \geq 2$. By (\ref{W11}) and Lemma~\ref{GN}, we get:
\begin{align} \nonumber
\lbrace \left\|u(t)\right\|_m^2 &\rbrace \overset{m}{\gtrsim} \lbrace \left\|u(t)\right\|_m^2  \left|u(t)\right|_{1,1}^{(4m-4)} \rbrace \overset{m}{\gtrsim} \lbrace \left\|u(t)\right\|^{4m-2}_1 \rbrace.
\end{align}
Thus, using H{\"o}lder's inequality and Lemma~\ref{finitetime}, we get:
\begin{align} \nonumber
\lbrace \left\|u(t)\right\|_m^2 &\rbrace \overset{m}{\gtrsim} \lbrace \left\|u(t)\right\|^{4m-2}_1 \rbrace \overset{m}{\gtrsim} \lbrace \left\|u(t)\right\|^{2}_1 \rbrace^{(2m-1)} \overset{m}{\gtrsim} \nu^{-(2m-1)}.\ \square
\end{align}
\\ \indent
The following two results generalise Lemma~\ref{finalexp}.

\begin{lemm} \label{finalexpbis}
For $m \geq 0$ and $p \in [1,\infty]$,
$$
\lbrace \left|u(t)\right|_{m,p}^2 \rbrace^{1/2} \overset{m,p}{\gtrsim} \nu^{-\gamma}.
$$
\end{lemm}

\textbf{Proof.}
The case $m=0$ is proved in Lemma~\ref{finitetimep}.
\\ \indent
In the case $m=1,\ p \geq 2$, it suffices to apply H{\"o}lder's inequality in place of Lemma~\ref{GN} in the proof of an analogue for Lemma~\ref{finalexp}.
\\ \indent
In the case $m \geq 2$, the proof is exactly the same as for Lemma~\ref{finalexp} for $p \in (1,\infty)$. In the cases $p=1,\infty$,  Lemma~\ref{GN} does not allow us to estimate $|u(t)|_{m,p}^2$ from below using $|u(t)|^2_{1,1}$ and $\Vert u(t) \Vert_1^2$. However, for $p=\infty$ we can proceed similarly, using the upper estimate (\ref{Lpupper}) for $|u(t)|^2_{\infty}$ and the lower estimate for $| u(t) |_{1,\infty}^2$. On the other hand, for $p=1$ it suffices to observe that we have $\left|u(t)\right|_{m,1} \geq \left|u(t)\right|_{m-1,\infty}$.
\\ \indent
Now consider the case $m=1,\ p \in [1,2)$. By H{\"o}lder's inequality we have:
\begin{align} \nonumber
\lbrace \left|u(t)\right|_{1,p}^2 \rbrace \geq & \lbrace \left\|u(t)\right\|_{1}^2  \rbrace^{2/p} \lbrace \left|u(t)\right|_{1,\infty}^2 \rbrace^{(p-2)/p}.
\end{align}
Using Lemma~\ref{finitetime} and Lemma~\ref{upperwmp}, we get the lemma's assertion. $\square$

\begin{lemm} \label{finalexpter}
For $m \geq 0$ and $p \in [1,\infty]$,
$$
\lbrace \left|u(t)\right|_{m,p}^{\alpha} \rbrace^{1/\alpha} \overset{m,p,\alpha}{\gtrsim} \nu^{-\gamma},\quad \alpha>0.
$$
\end{lemm}

\textbf{Proof.}
As previously, we may assume that $p>1$. The case $\alpha \geq 2$ follows immediately from Lemma~\ref{finalexpbis} and H{\"o}lder's inequality. The case $\alpha<2$ follows from H{\"o}lder's inequality, the case $\alpha=2$ and Lemma~\ref{upperwmp} (case $\alpha=3$), since we have:
\begin{align} \nonumber
\lbrace \left|u(t)\right|_{m,p}^{\alpha} \rbrace \geq & \lbrace \left|u(t)\right|_{m,p}^{2} \rbrace^{3-\alpha}
\lbrace \left|u(t)\right|_{m,p}^{3} \rbrace^{\alpha-2}.\ \square
\end{align}

The following theorem sums up the main results of this section, with the exception of Lemma~\ref{uxpos}.

\begin{theo} \label{avoir}
For $m \in \lbrace 0,1 \rbrace$ and $p \in [1,\infty]$, or for $m \geq 2$ and $p \in (1,\infty]$ we have:
\begin{equation} \label{avoirsim}
\Big( \lbrace \left|u(t)\right|_{m,p}^{\alpha} \rbrace \Big)^{1/\alpha} \overset{m,\alpha}{\sim} \nu^{-\gamma},\qquad \alpha>0,
\end{equation}
where $\lbrace \cdot \rbrace$ denotes time-averaging over $[T_1,T_2]$. The upper estimates in (\ref{avoirsim}) hold without time-averaging, uniformly for $t$ separated from $0$. Namely, we have:
$$
\left|u(t)\right|_{m,p} \overset{m,p}{\lesssim} \max(t^{-\gamma},\nu^{-\gamma}).
$$
On the other hand, the lower estimates hold for all $m \geq 0$ and $p \in [1,\infty]$. 
\end{theo}

\textbf{Proof.} 
Upper estimates follow from Lemma~\ref{upperwmp}, and lower estimates from Lemma~\ref{finalexpter}.\ $\square$

\section{Estimates for small-scale quantities} \label{turb}
\smallskip
\indent
In this section, we study analogues of quantities which are important for the study of hydrodynamical turbulence. We consider quantities in the physical space (structure functions) as well as in the Fourier space (energy spectrum). For notation for these quantities and the ranges $J_1,\ J_2,\ J_3$, see Subsection~\ref{notaturb}.
\\ \indent
Here, provided $\nu \leq \nu_0$, all estimates hold independently of the viscosity $\nu$. We recall that the brackets $\lbrace \cdot \rbrace$ stand for the averaging in time over an interval $[T_1,T_2]$: see (\ref{T1T2}).
\\ \indent
We begin by estimating the functions $S_{p}(\ell)$ from above.

\begin{lemm} \label{upperdiss}
For $\ell \in [0,1]$,
$$
S_{p}(\ell) \overset{p}{\lesssim} \left\lbrace \begin{aligned} & \ell^{p},\ 0 \leq p \leq 1. \\ & \ell^{p} \nu^{-(p-1)},\ p \geq 1. \end{aligned} \right.
$$
\end{lemm}

\textbf{Proof.} 
We begin by considering the case $p \geq 1$. We have:
\begin{align} \nonumber
S_{p}(\ell) &= \Big\{ \int_{S^1}{|u(x+\ell)-u(x)|^p dx} \Big\}
\\ \nonumber
& \leq \Big\{ \Big( \int_{S^1}{|u(x+\ell)-u(x)| dx} \Big) \Big( \max_{x} |u(x+\ell)-u(x)|^{p-1} \Big) \Big\}.
\end{align}
Using the fact that the space average of $u(x+\ell)-u(x)$ vanishes and H{\"o}lder's inequality, we obtain that
\begin{align} \nonumber
S_p(\ell) \leq & \Big\{ \Big(2 \int_{S^1}{(u(x+\ell)-u(x))^{+} dx} \Big)^{p} \Big\}^{1/p} \Big\{ \max_{x} |u(x+\ell)-u(x)|^{p} \Big\}^{(p-1)/p}
\\ \label{difference}
\leq & C \ell \Big\{ \max_{x} |u(x+\ell)-u(x)|^{p} \Big\}^{(p-1)/p},
\end{align}
where the second inequality follows from Lemma~\ref{uxpos}. Finally, by Theorem~\ref{avoir} we get:
\begin{align} \nonumber
S_p(\ell) & \leq C \ell \Big\{ ( \ell |u|_{1,\infty} )^{p} \Big\}^{(p-1)/p} \leq
C \ell^{p} \nu^{-(p-1)}.
\end{align}
The case $p<1$ follows immediately from the case $p=1$ since now $S_{p}(\ell) \leq (S_{1}(\ell))^p$, by H{\"o}lder's inequality. $\square$
\medskip
\\ \indent
For $\ell \in J_2 \cup J_3$, we have a better upper bound if $p \geq 1$.

\begin{lemm} \label{upperinert}
For $\ell \in J_2 \cup J_3$,
$$
S_{p}(\ell) \overset{p}{\lesssim} \left\lbrace \begin{aligned} & \ell^{p},\ 0 \leq p \leq 1. \\ & \ell,\ p \geq 1. \end{aligned} \right.
$$
\end{lemm}

\textbf{Proof.} The calculations are almost the same as in the previous lemma. The only difference is that we use another bound for the right-hand side of (\ref{difference}). Namely, by Theorem~\ref{avoir} we have:
\begin{align} \nonumber
S_{p}(\ell) & \leq C \ell \Big\{ \max_{x} |u(x+\ell)-u(x)|^{p} \Big\}^{(p-1)/p}
\\ \nonumber
& \leq C \ell \Big\{ (2 |u|_{\infty} )^{p} \Big\}^{(p-1)/p} \leq C \ell.\ \square
\end{align}

\begin{rmq}
The Lemmas~\ref{upperdiss} and \ref{upperinert} actually hold even if we drop the time-
\\
averaging, since in deriving them we only use upper estimates which hold uniformly for $t \geq T_1$.
\end{rmq}
\indent
To prove the lower estimates for $S_p(\ell)$, we need a lemma. Loosely speaking, this lemma states that there exists a large enough set $L_K \subset [T_1,T_2]$ such that for $t \in L_K$, several Sobolev norms are of the same order as their time averages. Thus, for $t \in L_K$, we can prove the existence of a \enquote{cliff} of height at least $C$ and width at least $C \nu$, using some of the arguments in \cite{AFLV92} which we exposed in the introduction.
\\ \indent
Note that in the following definition, (\ref{condi}-\ref{condii}) contain lower and upper estimates, while (\ref{condiii}) contains only an upper estimate. The inequality $|u(t)|_{\infty} \leq  \max u_x(t)$ in (\ref{condi}) always holds, since $u(t)$ has zero mean value and the length of $S^1$ is $1$.

\begin{defi}
For $K>1$, we denote by $L_K$ the set of all $t \in [T_1,T_2]$ such that the assumptions
\begin{align} \label{condi}
& K^{-1} \leq |u(t)|_{\infty} \leq  \max u_x(t) \leq K
\\ \label{condii}
& K^{-1} \nu^{-1} \leq  |u(t)|_{1,\infty} \leq K \nu^{-1}
\\ \label{condiii}
& |u(t)|_{2,\infty} \leq K \nu^{-2}
\end{align}
hold.
\end{defi}

\begin{lemm} \label{typical}
There exist constants $C,K_1>0$ such that for $K \geq K_1$, the
\\
Lebesgue measure of $L_K$ satisfies $\lambda(L_K) \geq C$.
\end{lemm}

\textbf{Proof.}
We begin by noting that if $K \leq K'$, then $L_K \subset L_{K'}$. By Lemma~\ref{uxpos} and Theorem~\ref{avoir}, for $K$ large enough the upper estimates in (\ref{condi}-\ref{condiii}) hold for all $t$. Therefore, if we denote by $B_K$ the set of $t$ such that
$$
\text{\enquote{The lower estimates in (\ref{condi}-\ref{condii}) hold for a given value of $K$}},
$$
then it suffices to prove the lemma's statement with $B_K$ in place of $L_K$. Now denote by $D_K$ the set of $t$ such that
$$
\text{\enquote{The lower estimate in (\ref{condii}) holds for a given value of $K$}}.
$$
By Lemma~\ref{GN} we have:
$$
|u|_{\infty} \geq C |u|_{2,\infty}^{-1} |u|_{1,\infty}^2.
$$
Thus if $D_K$ holds, then $B_{K'}$ holds for $K'$ large enough. Now it remains to show that there exists $C>0$ such that for $K$ large enough, we have the inequality $\lambda(D_K) \geq C$. We clearly have:
$$
\lbrace |u|_{1,\infty} \One(|u|_{1,\infty} < K^{-1} \nu^{-1}) \rbrace < K^{-1} \nu^{-1}.
$$
Here, $\One(A)$ denotes the indicator function of an event $A$. On the other hand, by the estimate for $\lbrace |u|_{1,\infty}^2 \rbrace$ in Theorem~\ref{avoir} we get:
\begin{align} \nonumber
\lbrace |u|_{1,\infty} \One( |u|_{1,\infty} > K \nu^{-1}) \rbrace & < K^{-1} \nu \lbrace |u|_{1,\infty}^2 \rbrace \leq C K^{-1} \nu^{-1}
\end{align}
Now denote by $f$ the function
$$
f=|u|_{1,\infty} \One(K_0^{-1} \nu^{-1} \leq |u|_{1,\infty} \leq K_0 \nu^{-1}).
$$
The inequalities above and the lower estimate for $\lbrace |u|_{1,\infty} \rbrace$ in Theorem~\ref{avoir} imply that
$$
\lbrace f \rbrace > (C-K_0^{-1}-C K_0^{-1}) \nu^{-1} \geq C_0 \nu^{-1},
$$
for some suitable constants $C_0$ and $K_0$. Since $f \leq K_0 \nu^{-1}$, we get:
$$
\lambda(f \geq C_0 \nu^{-1}/2) \geq C_0 K_0^{-1} (T_2-T_1)/2.
$$
Thus, since $|u|_{1,\infty} \geq f$, we have the inequality
$$
\lambda(|u|_{1,\infty} \geq C_0 \nu^{-1}/2) \geq C_0 K_0^{-1} (T_2-T_1)/2,
$$
which implies the existence of $C,K_1>0$ such that $\lambda(D_{K}) \geq C$ for $K \geq K_1$. $\square$
\smallskip
\\ \indent
Let us denote by $O_K \subset [T_1,T_2]$ the set defined as $L_K$, but with the relation (\ref{condii}) replaced by
\begin{equation} \label{condiibis}
K^{-1} \nu^{-1} \leq -\min u_x \leq K \nu^{-1}.
\end{equation}

\begin{cor} \label{typicalcor}
For $K \geq K_1$ and $\nu < K_1^{-2}$, we have $\lambda(O_K) \geq C$.
\end{cor}

\textbf{Proof.} For $K = K_1$ and $\nu < K_1^{-2}$, the estimates (\ref{condi}-\ref{condii}) tell us that
$$
\max u_x(t) \leq K_1 < K_1^{-1} \nu^{-1} \leq |u_x(t)|_{\infty},\quad t \in L_K.
$$
Thus, in this case we have $O_K=L_K$, which proves the corollary's assertion. Since increasing $K$ while keeping $\nu$ constant increases the measure of $O_K$, for $K \geq K_1$ and $\nu < K_1^{-2}$ we still have $\lambda(O_K) \geq C$. $\square$
\smallskip
\\ \indent
Now we fix
\begin{equation} \label{K}
K=K_1,
\end{equation}
and choose
\begin{equation} \label{nu0eq}
\nu_0=\frac{1}{6} K^{-2};\ C_1=\frac{1}{4}K^{-2};\ C_2=\frac{1}{20}K^{-4}.
\end{equation}
In particular, we have $0<C_1 \nu_0 <C_2<1$: thus the intervals $J_i$ are non-empty and non-intersecting for all $\nu \in (0,\nu_0]$. Everywhere below the constants depend on $K$.
\\ \indent
Actually, we can choose any values of $C_1$, $C_2$ and $\nu_0$, provided:
\begin{equation} \label{nu0ineq}
C_1 \leq \frac{1}{4}K^{-2};\quad 5 K^2 \leq \frac{C_1}{C_2}<\frac{1}{\nu_0}.
\end{equation}

\begin{lemm} \label{lowerdiss}
For $\ell \in J_1$,
$$
S_{p}(\ell) \overset{p}{\gtrsim} \left\lbrace \begin{aligned} &\ell^{p},\ 0 \leq p \leq 1. \\ & \ell^{p} \nu^{-(p-1)},\ p \geq 1. \end{aligned} \right.
$$
\end{lemm}

\textbf{Proof.}
By Corollary~\ref{typicalcor}, it suffices to prove that these upper estimates hold uniformly in $t$ for $t \in O_K$, with $S_{p}(\ell)$ replaced by
\begin{equation} \nonumber
\int_{S^1}{|u(x+\ell)-u(x)|^p dx}.
\end{equation}
Till the end of this proof, we assume that $t \in O_K$.
\\ \indent
Denote by $z$ the leftmost point on $S^1$ (considered as $[0,1)$) such that $ u'(z) \leq - K^{-1} \nu^{-1}$. Since $|u|_{2,\infty} \leq K \nu^{-2}$, we have
\begin{equation} \label{uxsmall}
u'(y) \leq -\frac{1}{2} K^{-1}  \nu^{-1},\quad y \in [z-\frac{1}{2} K^{-2} \nu,z+\frac{1}{2} K^{-2} \nu].
\end{equation}
In other words, the interval
$$
[z-\frac{1}{2} K^{-2} \nu,z+\frac{1}{2} K^{-2} \nu]
$$
corresponds to (a part of) a cliff.
\\ \indent
\textbf{Case $\mathbf{p \geq 1}$.} Since $\ell \leq C_1 \nu=\frac{1}{4}K^{-2} \nu$, by H{\"o}lder's inequality we get
\begin{align} \nonumber
\int_{S^1}&{|u(x+\ell)-u(x)|^p dx} \geq \int_{z-\frac{1}{4} K^{-2} \nu}^{z+\frac{1}{4} K^{-2} \nu}{|u(x+\ell)-u(x)|^p dx}
\\ \nonumber
&\geq (K^{-2} \nu/2)^{1-p} \Big( \int_{z-\frac{1}{4} K^{-2} \nu}^{z+\frac{1}{4} K^{-2} \nu}{|u(x+\ell)-u(x)| dx} \Big)^p
\\ \nonumber
&= C(p) \nu^{1-p} \Big( \int_{z-\frac{1}{4} K^{-2} \nu}^{z+\frac{1}{4} K^{-2} \nu}{ \Big(\int_{x}^{x+\ell}{- u'(y) dy } \Big) dx} \Big)^p
\\ \nonumber
&\geq C(p) \nu^{1-p} \Big( \int_{z-\frac{1}{4} K^{-2} \nu}^{z+\frac{1}{4} K^{-2} \nu}{ \frac{1}{2} \ell K^{-1} \nu^{-1} \ dx} \Big)^p = C(p) \nu^{1-p} \ell^p.
\end{align}
\\ \indent
\textbf{Case $\mathbf{p < 1}$.} By H{\"o}lder's inequality we obtain that
\begin{align} \nonumber
&\int_{S^1}{|u(x+\ell)-u(x)|^p dx} \geq \int_{S^1}{\Big((u(x+\ell)-u(x))^+\Big)^p dx}
\\ \nonumber
&\geq \Big( \int_{S^1}{\Big((u(x+\ell)-u(x))^+ \Big)^2 dx} \Big)^{p-1} \Big( \int_{S^1}{(u(x+\ell)-u(x))^+ dx} \Big)^{2-p}.
\end{align}
Using the upper estimate in (\ref{condi}) we get:
\begin{align} \nonumber
&\int_{S^1}{|u(x+\ell)-u(x)|^p dx}
\\ \nonumber
&\geq \Big( \int_{S^1}{\ell^2 K^2 dx} \Big)
^{p-1} \Big( \int_{S^1}{(u(x+\ell)-u(x))^+ dx} \Big)^{2-p}.
\end{align}
Since $\int_{S^1}{(u(\cdot+\ell)-u(\cdot))}=0$, we obtain that
\begin{align} \nonumber
&\int_{S^1}{|u(x+\ell)-u(x)|^p dx}
\\ \nonumber
&\geq C(p) \ell^{2 (p-1)} \Big(\frac{1}{2} \int_{S^1}{|u(x+\ell)-u(x)| dx} \Big)^{2-p} \geq C(p) \ell^p.
\end{align}
The last inequality follows from the case $p=1$.\ $\square$
\smallskip
\\ \indent
The proof of the following lemma uses an argument from \cite{AFLV92}, which becomes quantitative if we restrict ourselves to the set $O_K$.

\begin{lemm} \label{lowerinert}
For $m \geq 0$ and $\ell \in J_2$,
$$
S_{p}(\ell) \overset{p}{\gtrsim} \left\lbrace \begin{aligned} & \ell^{p},\ 0 \leq p \leq 1. \\ & \ell,\ p \geq 1. \end{aligned} \right.
$$
\end{lemm}

\textbf{Proof.} In the same way as above, it suffices to prove that the inequalities hold uniformly in $t$ for $t \in O_K$, with $S_{p}(\ell)$ replaced by
\begin{equation} \nonumber
\int_{S^1}{|u(x+\ell)-u(x)|^p dx},
\end{equation}
and we can restrict ourselves to the case $p \geq 1$. Again, till the end of this proof, we assume that $t \in O_K$.
\\ \indent
Define $z$ as in the proof of Lemma~\ref{lowerdiss}. We have
\begin{align} \nonumber
\int_{S^1}&{|u(x+\ell)-u(x)|^p dx} \geq
\\ \nonumber
&\int_{z-\frac{1}{2}\ell}^{z}
{ \Big| \underbrace{\int_{x}^{x+\ell}{u'^-(y)dy}}_{cliffs} - \underbrace{\int_{x}^{x+\ell}{u'^+(y)dy}}_{ramps} \Big|^p dx}.
\end{align}
Since $\ell \geq C_1 \nu=\frac{1}{4} K^{-2} \nu$, by (\ref{uxsmall}) for $x \in [z-\frac{1}{2}\ell,z]$ we get:
\begin{align} \nonumber
\int_{x}^{x+\ell}{u'^-(y)dy} &\geq \int_{z}^{z+\frac{1}{8} K^{-2} \nu}{u'^-(y) dy} \geq
\frac{1}{16} K^{-3}.
\\ \nonumber
&.
\end{align}
On the other hand, since $\ell \leq C_2$, by (\ref{condi}) and (\ref{nu0eq}) we get:
$$
\int_{x}^{x+\ell}{u'^+(y)dy} \leq C_2 K = \frac{1}{20} K^{-3}.
$$
Thus,
\begin{align} \nonumber
& \int_{S^1}{|u(x+\ell)-u(x)|^p dx} \geq \frac{1}{2} \ell \Bigg( \Big(\frac{1}{16}-\frac{1}{20}\Big) K^{-3} \Bigg)^p \geq C(p) \ell.\ \square
\end{align}
\medskip \\ \indent
Summing up the results above we obtain the following theorem.

\begin{theo} \label{avoir2}
For $\ell \in J_1$,
$$
S_{p}(\ell) \overset{p}{\sim} \left\lbrace \begin{aligned} & \ell^{p},\ 0 \leq p \leq 1. \\ & \ell^{p} \nu^{-(p-1)},\ p \geq 1. \end{aligned} \right.
$$
On the other hand, for $\ell \in J_2$,
$$
S_{p}(\ell) \overset{p}{\sim} \left\lbrace \begin{aligned} & \ell^{p},\ 0 \leq p \leq 1. \\ & \ell,\ p \geq 1. \end{aligned} \right.
$$
\end{theo}

The following result follows immediately from the definition (\ref{flatness}).

\begin{cor} \label{flatnesscor}
For $\ell \in J_2$, the flatness satisfies $F(\ell) \sim \ell^{-1}$.
\end{cor}

It remains to prove that, as long as $|k|$ remains in a certain range, after layer-averaging, we have  $\lbrace|\hat{u}(k)|^2 \rbrace \sim |k|^{-2}$. For this, we use a version of the Wiener-Khinchin theorem, stating that for any function $v \in L_2$ one has
\begin{equation} \label{spectrinertaux}
|v(\cdot+y)-v(\cdot)|^2=4\sum_{n \in \Z}{ \sin^2 (\pi ny) |\hat{v}(n)|^2}.
\end{equation}

\begin{theo} \label{spectrinert}
For $k$ such that $k^{-1} \in J_2$, we have $E(k) \sim k^{-2}$.
\end{theo}

\textbf{Proof.}
We recall that by definition (\ref{spectrum}),
$$
E(k) = \Bigg\{ \frac{\sum_{|n| \in [M^{-1}k,Mk]}{|\hat{u}(n)|^2}}{\sum_{|n| \in [M^{-1}k,Mk]}{1}} \Bigg\}.
$$
Therefore proving the assertion of the theorem is the same as proving that
\begin{equation} \label{spectrinertequiv}
\sum_{|n| \in [M^{-1}k,Mk]}{ n^2 \lbrace |\hat{u}(n)|^2 \rbrace } \sim k.
\end{equation}
From now on, we will indicate explicitly the dependence on $M$. The upper estimate holds without averaging over $n$ such that $|n| \in [M^{-1} k, Mk]$. Indeed,  by (\ref{intparts}) we know that
\begin{equation} \nonumber
\lbrace |\hat{u}(n)|^2 \rbrace \leq C n^{-2}.
\end{equation}
Also, this inequality implies that
\begin{equation} \label{spectrinertupper1}
\sum_{|n| < M^{-1} k}{ n^2 \lbrace |\hat{u}(n)|^2 \rbrace } \leq C M^{-1} k
\end{equation}
and
\begin{equation} \label{spectrinertupper2}
\sum_{|n| > M k}{\lbrace |\hat{u}(n)|^2 \rbrace } \leq C M^{-1} k^{-1}.
\end{equation}
Now it remains to prove the lower bound. We have:
\begin{align} \nonumber
\sum_{|n| \leq M k}{ n^2 \lbrace |\hat{u}(n)|^2 \rbrace } & \geq \frac{k^2}{\pi^2} \sum_{|n| \leq M k}{ \sin^2 (\pi nk^{-1}) \lbrace |\hat{u}(n)|^2 \rbrace }
\\ \nonumber
&\geq \frac{k^2}{\pi^2} \Big( \sum_{n \in \Z}{ \sin^2 (\pi nk^{-1}) \lbrace |\hat{u}(n)|^2 \rbrace } - \sum_{|n| > M k}{\lbrace |\hat{u}(n)|^2 \rbrace } \Big).
\end{align}
Using (\ref{spectrinertaux}) and (\ref{spectrinertupper2}) we get:
\begin{align} \nonumber
\sum_{|n| \leq M k}{ n^2 \lbrace |\hat{u}(n)|^2 \rbrace } &\geq  \frac{k^2}{4 \pi^2} \Big( \lbrace |u(\cdot+k^{-1})-u(\cdot)|^2 \rbrace - C M^{-1} k^{-1} \Big)
\\ \nonumber
&\geq \frac{k^2}{4 \pi^2} \Big( S_2(k^{-1})-C M^{-1} k^{-1} \Big).
\end{align}
Finally, using Theorem~\ref{avoir2} we obtain that
\begin{equation} \nonumber
\sum_{|n| \leq M k}{ n^2 \lbrace |\hat{u}(n)|^2 \rbrace }  \geq (C-C M^{-1}) k.
\end{equation}
Now we use (\ref{spectrinertupper1}) and we choose $M \geq 1$ large enough to obtain (\ref{spectrinertequiv}). $\square$

\section{Estimates for small-scale quantities in the inviscid limit} \label{inviscid}
\smallskip
 \indent
It is a well-known fact (see for instance \cite{KrLo89}) that as $\nu$ tends to $0$, the solutions of (\ref{Burbegin}) converge to \textit{weak entropy solutions} of the inviscid equation $u_t+f'(u)u_x=0$, for fixed $t$. The convergence takes place for almost every $x$, and therefore also in $L_1$, since solutions are uniformly bounded for all $\nu$.
\\ \indent
These solutions, denoted $u^0$, inherit all previously proved properties which hold uniformly for small enough $\nu$ for solutions of (\ref{Burbegin}).
\\ \indent
To begin with, we define the non-empty and non-intersecting intervals
\begin{equation} \nonumber
J_2=(0,\ C_2];\ J_3=(C_2,\ 1],
\end{equation}
which now correspond to the inertial range and the energy range, respectively. The inviscid Burgers equation does not have a dissipation range, since formally there is no dissipation, despite the presence of an anomaly due to the shocks \cite{BK07}. The constant $C_2$ is the same as above.
\\ \indent
Then we define $S_p^0$, $F^0$ and $E^0$ for solutions $u^0(t,x)$ in the same way as the previously considered quantities $S_p$, $F$ and $E$ for solutions of the viscous equation. By the dominated convergence theorem, we obtain the following results:

\begin{theo} \label{avoir2inv}
For $\ell \in J_2$, $S^0_{p}(\ell) \overset{p}{\sim} \left\lbrace \begin{aligned} & \ell^{p},\ 0 \leq p \leq 1. \\ & \ell,\ p \geq 1. \end{aligned} \right.$
\end{theo}

\begin{cor} \label{flatnesscorinv}
For $\ell \in J_2$, the flatness satisfies $F^0(\ell) \sim \ell^{-1}$.
\end{cor}

\begin{theo} \label{spectrinertinv}
For $k$ such that $k^{-1} \in J_2$, we have $E^0(k) \sim k^{-2}$.
\end{theo}

\section*{Acknowledgements}
\smallskip
 \indent
I would like to thank my Ph.D. advisor S.Kuksin, who originally formulated the problem. I am also very grateful to A.Biryuk, in particular for bringing to my attention some articles of S.Kruzhkov, and to J.Bec, U.Frisch, F.Golse and K.Khanin for many helpful discussions. Parts of the present work were done during my stays at Laboratoire AGM, University of Cergy-Pontoise, INRIA-IRISA-ENS Cachan Bretagne and DPT, University of Geneva. I have been supported respectively by the grants ERC BLOWDISOL, GEOPARDI and BRIDGES: I would like to thank all the faculty and staff, and especially the principal investigators, respectively F.Merle, E.Faou and J.-P.Eckmann, for their hospitality.
\bigskip

\bibliographystyle{plain}
\bibliography{Bibliogeneral}

\end{document}